\documentclass[aos,MSNbibl,citesort,dvips]{arximspdf}
\usepackage{graphicx}

\doi{10.1214/12-AOS984} 
\referstodoi{10.1214/11-AOS949}
\volume{40}
\issue{4}
\pubyear{2012}
\firstpage{1984}
\lastpage{1988}

\makeatletter
\newcommand{\argmin}{\operatorname{argmin}}
\newcommand{\X}{\mathbf{X}}
\newcommand{\hS}{\widehat{S}}
\newcommand{\hL}{\widehat{L}}
\makeatother

\begin{document}
\begin{frontmatter}

\title{Discussion: Latent variable graphical model selection via convex optimization}
\runtitle{Comment}

\begin{aug}
\author[A]{\fnms{Christophe} \snm{Giraud}\corref{}\ead[label=e1]{Christophe.Giraud@polytechnique.edu}}
\and
\author[B]{\fnms{Alexandre} \snm{Tsybakov}\ead[label=e2]{Alexandre.Tsybakov@ensae.fr}}
\runauthor{C. Giraud and A. Tsybakov}
\affiliation{Ecole Polytechnique and CREST-ENSAE}
\address[A]{CMAP, UMR CNRS 7641\\
Ecole Polytechnique\\
Route de Saclay \\
F-91128 Palaiseau Cedex\\
France\\
\printead{e1}} 
\address[B]{Laboratoire de Statistique\\
CREST-ENSAE \\
3, av. Pierre Larousse \\
F-92240 Malakoff Cedex\\
France \\
\printead{e2}}
\end{aug}

\received{\smonth{2} \syear{2012}}



\end{frontmatter}

Recently there has been an increasing interest in the problem of
estimating a high-dimensional matrix $K$ that can be decomposed in a
sum of a sparse matrix $S^*$ (i.e., a matrix having only a small number
of nonzero entries) and
a low rank matrix $L^*$. This is motivated by applications in computer
vision, video segmentation, computational biology, semantic indexing,
etc. 
The main contribution and novelty of the Chandrasekaran, Parrilo and
Willsky paper (CPW in what follows) is to propose and study a method of
inference about such decomposable matrices for a particular setting
where $K$ is the precision (concentration) matrix of a partially
observed sparse Gaussian graphical model (GGM). In this case, $K$ is
the inverse of the covariance matrix of a Gaussian vector $X_{O}$
extracted from a larger Gaussian vector $(X_{O},X_{H})$ with sparse
inverse covariance matrix. Then it is easy to see that $K$ can be
represented as a sum of a sparse precision matrix $S^*$ corresponding
to the observed variables $X_{O}$ and a matrix $L^*$ with rank at most
$h$, where $h$ is the dimension of the latent variables $X_{H}$. If $h$
is small, which is a typical situation in practice, then $L^*$ has low
rank. The GGM with latent variables is of major interest for
applications in biology or in social networks where one often does not
observe all the variables relevant for depicting sparsely the
conditional dependencies. Note that formally this is just one possible
motivation and mathematically the problem is dealt with in more
generality, namely, postulating that the precision matrix satisfies
%
%
\begin{equation}\label{1}
K=S^*+L^*
\end{equation}
with sparse $S^*$ and low-rank $L^*$, both symmetric matrices. A small
amendment to that inherited from the latent variables motivation is
that $L^*$ is assumed negative definite (in our notation, $L^*$
corresponds to $-L^*$ in the paper). We believe that this is not
crucial and all the results remain valid without this assumption.

CPW propose\vspace*{1pt} to estimate the pair $(S^*,L^*)$ from a $n$-sample of
$X_{O}$ by the pair $(\widehat{S},\widehat{L})$ obtained by
minimizing the negative log-likelihood with mixed $\ell^1$ and nuclear
norm penalties; cf. (1.2) of the paper. The key issue in this context
is identifiability. Under what conditions can we identify $S^*$ and
$L^*$ separately? CPW provide geometric conditions of identifiability
based on transversality of tangent spaces to the varieties of sparse
and low-rank matrices. They show that, under these conditions, with
probability close to 1, it is possible to recover the support of $S^*$,
the rank of $L^*$ and to get a bound of order $O(\sqrt{p/n})$ on the
estimation errors $|\widehat{S}- S^*|_{\ell^\infty}$ and
$\|\widehat{L}- L^*\|_2$. Here, $p$ is the dimension of $X_{O}$ and
$|\cdot|_{\ell^q}$ and $\|\cdot\|_2$ stand for the componentwise
$\ell^q$-norm and the spectral norm of a matrix, respectively.

Overall, CPW pioneer a hard and important problem of high-dimensional
statistics and provide an original solution both in the theory and in
numerically implementable realization. While being the first work to
shed light on the problem, the paper does not completely raise the
curtain and several aspects still remain to be understood and
elucidated.

\section*{The nature of the results} The most important problem for current
applications appears to be the estimation of $S^*$ or the recovery of
its support. Indeed, the main interest is in the conditional
dependencies of the coordinates of $X_{O}$ in the complete model
$(X_{O},X_{H})$ and this information is carried by the matrix $S^*$. In
this context, $L^*$ is essentially a nuisance, so that bounds on the
estimation error of $L^*$ and the recovery of the rank of $L^*$ are of
relatively moderate interest. However, mathematically, the most
sacrifice comes from the desire to have precise estimates of $L^*$.
Indeed, if $\widehat\Sigma_n$ and $\Sigma$ denote the empirical and
population covariance matrices,
the slow rate $O(\sqrt{p/n})$ comes from the bound on $\|\widehat
\Sigma_n-\Sigma\|_{2}$ in Lemma 5.4, that is, from the stochastic
error corresponding to $L^*$. Since the sup-norm error $|\widehat
\Sigma_n-\Sigma|_{\ell^\infty}$ is of order $\sqrt{(\log p)/n}$,
can we get a better rate when solely focusing on $|\hS-S^*|_{\ell
^{\infty}}$?

\section*{Extension to high dimensions} The results of the paper are valid
and meaningful only when $p<n$. However, for the applications of GGM,
the case $p\gg n$ is the most common. A key question is whether the
restriction $p<n$ is intrinsic, that is, whether it is possible to have
results on $S^*$ in model  (\ref{1}) when $p\gg n$. Since the
traditional model with sparse component $S^*$ alone is still tractable
when $p\gg n$, a related question is whether introducing the model
(\ref{1}) with two components and estimating both $S^*$ and $L^*$
gives any improvement in the $p\gg n$ setting as compared to
estimation in the model with a sparse component alone.
A small simulation study that we provide below suggests that already
for $p=n,$ including the low-rank component \textit{in the estimator} may
yield no improvement as compared to traditional sparse estimation
without the low-rank component, although this low-rank component is
effectively present \textit{in the model}.

\section*{Optimal rates} The paper obtains bounds of order $O(\sqrt{p/n})$
on the estimation errors $|\widehat{S}- S^*|_{\ell^\infty}$ and
$\|\widehat{L}- L^*\|_2$ with probability $1-2\exp(-p)$.
Can we achieve a better rate than $\sqrt{p/n}$ when solely focusing on
the recovery of $S^*$ with the usual probability $1-p^{-a}$ for some
$a>0$? Is the rate $\sqrt{p/n}$ optimal in a minimax sense on some
class of matrices? Note that one should be careful
in defining the class of matrices because in reality the rate is not
$O(\sqrt{p/n})$ but rather $O(\psi\sqrt{p/n})$, where $\psi$ is the
spectral norm of $\Sigma$ depending on $p$. It can be large for large
$p$. Surprisingly, not much is known about the optimal rates even in
the simpler case of purely sparse precision matrices, without the
low-rank component. In this case, \cite{RavikumarEtal,cai_etal} and
\cite{sun_zhang} provide some analysis of the upper
bounds on the estimation error of different estimators and under
different sets of assumptions on the precision matrix. All these bounds
are of ``order'' $O(\sqrt{(\log p)/n})$, but again one should be very
careful here because of the factors depending on $p$ that multiply this
rate. In \cite{cai_etal}, the factor is the squared $\ell^1 \to\ell
^1$ norm of the precision matrix while in \cite{RavikumarEtal}, it is
the squared degree of the graphical model multiplied by some
combinations of powers of matrix norms that are not easy to interpret.
The most recent paper \cite{sun_zhang} obtains the rate $O(d\sqrt
{(\log p)/n})$, where $d$ is the degree of the graph for $\ell^\infty
$-bounded precision matrices. An open problem is to find optimal rates
of convergence
on classes of precision matrices defined via sparsity and low rank
characteristics. The same problem makes sense for covariance matrices.
Here, some advances have been achieved very recently. In particular,
some optimal rates of estimation of low-rank covariance matrices are
provided by \cite{Lounici}.

\textit{The assumptions} of the paper are stated in terms of some
inaccessible characteristics such as $\xi(T)$ and $\mu(\Omega)$ and
seem to be very strong. They are in the spirit of the
irrepresentability condition for the vector case used to prove model
selection consistency of the Lasso. For a given set of data, there is
no means to check whether these assumptions are satisfied. What happens
when they do not hold? Can we still have some convergence properties
under no assumption at all or under weaker assumptions akin to the
restricted eigenvalue condition in the vector case?

\section*{Choice of the tuning parameters} The choice of parameters
$(\gamma,\lambda_{n})$ ensuring algebraic consistency in Theorem 4.1
depends on various unknown quantities. Proposing a reasonable
data-driven selector for $(\gamma,\lambda_{n})$ (e.g., similarly
to \cite{GHV12} for the pure sparse setting) would be very helpful for
the practice.

\section*{Alternative methods of estimation} Constructively, the method of
CPW is obtained from the GLasso of \cite{FHT} by adding a penalization
by the nuclear norm of the low-rank component. Similar low-rank
extensions can be readily derived from other methods, such as the
Dantzig type approach of \cite{cai_etal} and the regression approach of
\cite{MB06,Gir08}. Consider a Gaussian random vector $X\in\mathbb{R}^p$
with mean 0 and nonsingular covariance matrix $\Sigma$. Let
$K=\Sigma^{-1}$ be the precision matrix. We assume that $K$ is of the
form (\ref{1}) where $S^*$ is sparse and $L^*$ has low rank.

(a) \textit{Dantzig type approach.} In the spirit of \cite{cai_etal}, we
may define our estimator as a solution of the following convex program:
%
%
\begin{equation}\label{critSL1}
(\hS,\hL)=\argmin\limits_{(S,L)\in\mathcal{G}} \{ |S|_{\ell^1}+\mu\|L\|_{*}
\},
\end{equation}
where $\|\cdot\|_{*}$ is the nuclear norm, $\mathcal{G} = \{(S,L)\dvtx
|\widehat\Sigma_{n}(S+L)-I|_{\ell^\infty} \le\lambda\}$ and $\mu,
\lambda>0$ are tuning constants. Here, the nuclear norm $\|L\|_{*}$ is
a convex relaxation of the rank of $L^*$.

(b) \textit{Regression approach.} The regression
approach \cite{MB06,Gir08} is an alternative point of view for
estimating the structure of a GGM. In the pure sparse setting, some
numerical experiments \cite{VillersEtal} suggest that it may be more
reliable than the $\ell^1$-penalized log-likelihood approach.
Let $\operatorname{diag}(A)$ denote the diagonal of square matrix $A$ and $\|
A\|_F$ its Frobenius norm. Defining
\[
\Theta=\argmin\limits_{A : \operatorname{diag}(A)=0}\|\Sigma^{1/2}(I-A)\|_F^2,
\]
we have $\Theta=K\Delta+I$, where $I$ is the identity matrix and
$\Delta$ is the diagonal matrix with diagonal elements
$\Delta_{jj}=-1/K_{jj}$ for $j=1,\ldots,p$. Thus, we have the
decomposition
\[
\Theta=\bar{S}+\bar{L},\qquad \mbox{where } \bar{S}=S^*\Delta+I\
\mbox{ and } \bar{L}=L^*\Delta.
\]
Note that $\operatorname{rank}(\bar{L})=\operatorname{rank}(L^*)$ and the nondiagonal
elements $\bar{S}_{ij}$ of matrix $\bar{S}$ are nonzero only if
$S^*_{ij}$ is nonzero. Therefore, recovering the support of $S^*$ and
$\operatorname{rank}(L^*)$ is equivalent to recovering the support of $\bar{S}$
and $\operatorname{rank}(\bar{L})$.

%
\begin{figure}

\includegraphics{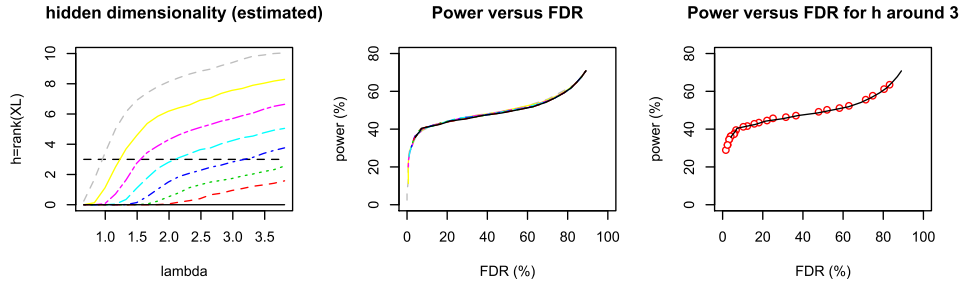}

\caption{Each color corresponds to a fixed value of $\mu$, the
solid-black color being for $\mu=+\infty$. For each choice of $\mu$,
different quantities are plotted for a series of values of~$\lambda$.
Left: Mean rank of $\X\hL$. Middle: The curve of estimated power
versus estimated FDR. Right: The power versus FDR for the estimators
fulfilling $\mathbb{E}[\operatorname{rank}(\X\hL)]\approx h=3$ (red dots),
superposed with the Power versus the FDR for $\mu=+\infty$ (in
solid-black).}\label{fig1}
\end{figure}

Now, we estimate $(\bar{S},\bar{L})$ from an $n$-sample of $X$
represented as an $n\times p$ matrix~$\X$. Noticing that the sample
analog of $\|\Sigma^{1/2}(I-A)\|_F^2$ is $\|\X(I-A)\|_F^2/n$ and using
the decomposition $\Theta=\bar{S}+\bar{L}$, we arrive at the
following estimator:
%
%
\begin{equation}\label{critSL}
\hspace*{10pt}(\hS,\hL)=\argmin\limits_{(S,L): \operatorname{diag}(S+L)=0}\biggl\{\frac{1}{2}\|\X
(I-S-L)\|_F^2+\lambda|S|_{\ell^1, \mathrm{off}}+\mu\|\X L\|_{*}\biggr\},
\end{equation}
where $\mu, \lambda$ are positive tuning constants and $|S|_{\ell^1,
\mathrm{off}}=\sum_{i\neq j}|S_{ij}|$. Note that here the low-rank
shrinkage is driven by the nuclear norm $\|\X L\|_{*}$ rather than by
$\|L\|_{*}$.
The convex minimization in (\ref{critSL}) can be performed efficiently
by alternating block descents on the off-diagonal elements of $S$, the
matrix $L$ and the diagonal of $S$. The off-diagonal support of $S^*$
is finally estimated by the off-diagonal support of~$\hS$.

\section*{Numerical experiment} A sparse Gaussian graphical model in
$\mathbb{R}^{30}$ is generated randomly according to the procedure
described in Section 4 of \cite{GHV12}. A~sample of size $n=30$ is
drawn from this distribution and $\X$ is obtained by hiding the values
of 3 variables. These 3 hidden variables are chosen randomly among the
connected variables. The estimators $(\hS,\hL)$ defined
in~(\ref{critSL}) are then computed for a grid of values of $\lambda$ and
$\mu$. The results are summarized in Figure \ref{fig1} (average over 100
simulations).

Strikingly, there is no significative difference in these examples
between the procedure of 
\cite{MB06} (corresponding to $\mu=+\infty$, in solid-black) and the
procedure (\ref{critSL}) that includes the low-rank component
(corresponding to finite~$\mu$).


%

\printaddresses

\end{document}